\documentclass[12pt]{amsart}

\usepackage{amsmath}
\usepackage{amsthm}
\usepackage{amsfonts}
\usepackage{amssymb,amsmath,latexsym}
\usepackage{graphicx}

\theoremstyle{plain}
\newtheorem{theorem}{Theorem}[section]
\newtheorem{lemma}[theorem]{Lemma}
\newtheorem{proposition}[theorem]{Proposition}
\newtheorem{corollary}[theorem]{Corollary}

\theoremstyle{definition}
\newtheorem{remark}[theorem]{Remark}

\title[SLE(8/3) in annuli]{Restricting SLE(8/3) to an annulus}

\author{Robert~O. Bauer}
\address{Altgeld Hall\\Department of Mathematics\\ 
	University of Illinois at Urbana-Champaign\\ 
	1409 West Green Street \\ 
	Urbana, IL 61801, USA}
\email{rbauer@math.uiuc.edu}

\begin{document}

\begin{abstract}
We study the probability that chordal $\text{SLE}_{8/3}$ in the unit disk from $\exp(ix)$ to $1$ avoids the disk of radius $q$ centered at zero. We find the initial/boundary-value problem satisfied by this probability as a function of $x$ and $a=\ln q$, and show that asymptotically as $q$ tends to one this probability decays like $\exp(-cx/(1-q))$ with $c=5\pi/8$ for $0< x\le\pi$. We also give a representation of this probability as a multiplicative functional of a Legendre process.   
\end{abstract}

\maketitle

\section{Introduction}
In this paper we study certain hitting probabilities for the chordal Schramm-Loewner evolution with parameter $\kappa=8/3$ ($\text{SLE}_{8/3}$). We study this question for $\text{SLE}_{8/3}$ because this process lies in the intersection of two important classes of conformally invariant measures. 

On the one hand, we have chordal SLE: these are families of measures on non-self-crossing curves $\gamma$, indexed by the simply connected domain $D$ the curve $\gamma$ lives in, and the endpoints $z,w$ of $\gamma$ on $\partial D$. We can think of $\gamma$ as a random interface separating two different materials on $D$. If $P_{D,z\to w}$ denotes the law of the curve $\gamma$ in $D$ from $z$ to $w$, then the family $\{P_{D,z\to w}\}$ is a Schramm-Loewner evolution if members of the family are related by

(1) {\em conformal invariance:} if $f$ is a conformal map from $D$ to $ D'$ sending $z,w$ to $z',w'$, then $f\circ P_{D,z\to w}=P_{D',z'\to w'}$;

(2) {\em domain Markovianity:} if  $\gamma$ has law $P_{D,z\to w}$, $z'$ is an interior point of $\gamma$, and we condition on the  segment $\gamma'$ of $\gamma$ from $z$ to $z'$, then the remaining segment of $\gamma$, from $z'$ to $w$, has law $P_{D\backslash\gamma',z'\to w}$;

and if, for the particular case where $D$ is the upper half-plane $\mathbb H$, $z=0$, $w=\infty$, the law of $\gamma$ is symmetric with respect to the imaginary axis. Suppose $\{P_{D,z\to w}\}$ is  such a family. Using L\"owner's theory of slit mappings \cite{loewner}, Schramm showed that if $t\in[0,\infty)\mapsto\gamma_t\in\overline{\mathbb H}$ is correctly parameterized, $\gamma_0=0$, $D_t$ denotes the unbounded component of $\mathbb H\backslash\gamma(0,t]$, and $g_t:D_t\to\mathbb H$ is conformal with `hydrodynamic' normalization at infinity \[\lim_{z\to\infty}g_t(z)-z=0,\] then, under $P_{\mathbb H,0\to\infty}$, $g_t(\gamma_t)=\sqrt{\kappa}B_t$ for a standard 1-dimensional Brownian motion $\{B_t:t\ge0\}$ starting at zero and a constant $\kappa\ge0$, \cite{schramm:SLE}.

On the other hand, we have {\em restriction measures}. These are again families of measures $\{P_{D,z,w}\}$ indexed by simply connected domains $D$ and two boundary points $z,w$, but this time describing random, closed, simply connected subsets (which we denote also by $\gamma$) of $\overline{D}$ such that $\gamma\cap\partial D=\{z,w\}$. For example, a simple curve in $D$ from $z$ to $w$ is such a set. We dropped the $\to$ in the notation as $\gamma$ is a point-set without a `direction.' A family $\{P_{D,z,w}\}$ is called a restriction measure if it is conformally invariant (as in (1) above), and satisfies the 

(3) {\em restriction property:} if $\gamma$ has distribution $P_{D,z, w}$, $D'\subset D$ and $z,w\in\partial D'$, then conditional on $\{\gamma\subset D'\}$ the distribution of $\gamma$ is $P_{D',z,w}$. 

In the statement of the restriction property it is understood that $z$ and $w$ are bounded away from the part of the boundary of $D$ that does not belong to $\partial D'$. An example of a restriction measure is provided by the `filling' of a Brownian excursion in $D$ from $z$ to $w$. Restriction is a powerful property. If $\{P_{D,z,w}\}$ denotes a restriction measure,  and if $D_2\subset D_1\subset D$ and $z,w\in\partial D_2$, then restriction implies in particular that
\begin{equation}\label{E:rest1}
P_{D,z, w}\{\gamma\subset D_2\}=P_{D_1,z, w}\{\gamma\subset D_2\}P_{D,z, w}\{\gamma\subset D_1\}.
\end{equation}
By conformal invariance it is enough to consider the case when $D$ is the upper half-plane $\mathbb H$, $z=0$, and $w=\infty$. That is, suppose that $D_2\subset D_1\subset\mathbb H$ with $0,\infty\in\partial D_2$. Denote $\Phi_{1,2}:D_{1,2}\to\mathbb H$ the conformal map with normalization $\lim_{z\to\infty}\Phi_{1,2}(z)/z=1$, $\Phi_{1,2}(0)=0$. Then we can rewrite \eqref{E:rest1} as
\begin{equation}\label{E:rest2}
P_{\mathbb H,0,\infty}\{\gamma\subset D_2\}=P_{\mathbb H,0,\infty}\{\gamma\subset\Phi_1(D_2)\}P_{\mathbb H,0,\infty}\{\gamma\subset D_1\}
\end{equation}
As we can identify a domain with the unique normalized conformal map from that domain to $\mathbb H$, we may write $F(\Phi_{1,2})=P_{\mathbb H,0,\infty}\{\gamma\subset D_{1,2}\}$. In particular,  \eqref{E:rest2} is equivalent to
\begin{equation}\label{E:rest3}
F(\Phi_2)=F(\Phi_2\circ\Phi_1^{-1})\cdot F(\Phi_1),
\end{equation}
that is, $F$ is a homomorphism from the semigroup of conformal maps (with composition) to $[0,\infty)$ (with multiplication). Lawler, Schramm, and Werner showed that this implies the remarkable result that there exists an $\alpha>0$ such that 
\begin{equation}\label{E:rest4}
P_{\mathbb H,0,\infty}\{\gamma\in D\}=F(\Phi)=\Phi'(0)^{\alpha},
\end{equation}
where $D$ is a {\em simply connected\ } subdomain of $\mathbb H$ containing $0,\infty$ as boundary points, see \cite{LSW:rest}. If $\gamma$ is both, an SLE and a restriction measure, then 
\begin{align}\label{E:rest5}
P_{\mathbb H,0\to\infty}\{\gamma\subset D|\gamma[0,t]\}&=1\{\gamma[0,t]\subset D\}P_{\mathbb H,W_t\to\infty}\{\gamma\subset g_t(D)\}\notag\\
&=1\{\gamma[0,t]\subset D\}h_t'(W_t)^{\alpha},
\end{align}
where $h_t$ is the normalized conformal map from $g_t(D)$ to $\mathbb H$, and $W_t=\sqrt{\kappa}B_t$. The first equality in \eqref{E:rest5} is on account of $\gamma$ being an SLE, the second a consequence of restriction. It follows that $h_t'(W_t)^{\alpha}$ is a martingale on $\{\gamma[0,t]\subset D\}$. A calculation now shows that this implies $\kappa=8/3$ and $\alpha=5/8$, \cite{LSW:rest}. The self-avoiding random walk satisfies the discrete version of the restriction property and it is conjectured that the scaling limit of self-avoiding random walk is $\text{SLE}_{8/3}$, \cite{lsw:saw}.  

We now ask what happens if we restrict to domains $D\subset\mathbb H$ with `holes,' i.e if $D$ is no longer simply connected. Then there is no homeomorphism  from $D$ to $\mathbb H$. Even more, while connectivity classifies topological equivalence, it does not classify conformal equivalence. For example, two annuli are conformally equivalent if and only if the ratio of outer to inner radius of the former equals that of the latter. In other words, there is a conformal parameter, or modulus, which labels conformal equivalence classes of doubly connected domains, \cite{ahlfors:1978}.

However, it is easy to extend restriction measures to multiply connected domains. Suppose $\{P_{D,zw}\}$ is a restriction measure as above. If $D'$ is finitely connected and $z,w$ points on the same boundary component of $D'$, we define 
\begin{equation}\label{E:beffara}
P_{D',z,w}=P_{D,z,w}\{\ \cdot\ |\gamma\subset D'\},
\end{equation}
where $D\supset D'$ is simply connected and $z,w\in\partial D$. Restriction for simply connected domains implies that $P_{D',z,w}$ is independent of the choice of $D$, and an inclusion/exclusion argument of Beffara shows that then \eqref{E:beffara} holds for arbitrary finitely connected domains $D',D$ with $D'\subset D$, $z,w\in\partial D'\cap\partial D$, \cite{beffara}. The identity \eqref{E:rest1} still holds in this more general context but \eqref{E:rest2} and \eqref{E:rest4} no longer make sense. Thus, while we can define restriction measures in multiply connected domains, we cannot calculate---or do not have a functional expression for---the probability that $\gamma$ hits a `hole.' Finding a functional expression which generalizes \eqref{E:rest4} to multiply connected domains is the main motivation for this paper.  

To begin, we decided to focus on the simplest case, just one hole, and address this case for the restriction measure which also is an SLE, making SLE-tools available. So suppose $\gamma$ is a chordal $\text{SLE}_{8/3}$ in the unit disk $\mathbb U=\{|z|<1\}$ from $e^{ix}$ to $1$ and $A_q=\{q<|z|<1\}$ an annulus. Then
\[
	P_{\mathbb U,e^{ix}\to1}\{\gamma\subset A_q\}
\]
is a function $F$ of $x$ and $a=\ln q$. In this paper we show that $F$ is $C^{1,2}$, find the initial/boundary-value problem to which this function is the solution, see Theorem \ref{T:main}, and show in Theorem \ref{T:asymp0} that asymptotically
\begin{equation}\label{E:asymptotics}
F(a,x)\asymp\exp\left(-\frac{5\pi}{8}\cdot\frac{x}{1-q}\right),\quad 0\le x\le \pi,
\end{equation}
as $q\nearrow1$. Using this strong decay we obtain a stochastic representation for $F(a,x)$  as
\begin{align}
&\left[\prod_{n=1}^{\infty}\frac{1-2q^{2n}+q^{4n}}{1-2q^{2n}\cos x+q^{4n}}\right]^{3/4}\notag\\
&\quad\times\mathbb E\left[\exp\left(\int_a^{\sigma}\left[\frac{1}{12}-\sum_{n=1}^{\infty}\frac{2n e^{2nb}}{1-e^{2nb}}\left(1-\cos n Y_b\right)\right] db\right),\sigma<0\right]
\end{align}
in Theorem \ref{T:rep}. Here $Y$ is a Legendre process on $[0,2\pi]$ starting at $x$ at time $a<0$ and $\sigma$ is the first time $Y$ hits the boundary. We give an alternative expression in terms of Jacobi's $\vartheta$-function and Weierstrass' $\wp$-function. 

In \cite{werner:2005}, Werner also studies the asymptotics of a non-intersection probability in annuli as $q\nearrow 1$, namely the probability, appropriately rescaled, that chordal $\text{SLE}_{8/3}$ from ``near $1$'' to $1$ in the unit disk stays in the annulus $A_q$ and goes the long way (around the hole), see \cite[Lemma 18]{werner:2005}. He finds that that probability decays like $\exp(-5\pi^2/(4(1-q))$. This result can be guessed from \eqref{E:asymptotics} as follows. The probability that a chordal $\text{SLE}_{8/3}$ from ``near 1'' to $1$ goes around the disk of radius $q$ centered at zero is, for $q$ close to 1, approximately the same as the probability that a chordal $\text{SLE}_{8/3}$ from $1$ to $-1$ goes around the disk of radius $q$ via the upper half-plane, followed by an independent $\text{SLE}_{8/3}$ from $-1$ to $1$, which goes around the disk of radius $q$ via the lower half-plane. Thus the probability Werner calculates should behave asymptotically like the square of \eqref{E:asymptotics} for $x=\pi$, which indeed is the case.

Concerning the behavior of $F(a,x)$ as $q\searrow0$ a brief analysis of the initial/boundary-value problem leads to the conjecture
\begin{equation}
	F(a,x)=1-cq^{2/3}\sin^2 x/2,\quad q\searrow 0,
\end{equation}
for some constant $c$, see Proposition \ref{P:asymp}. We give evidence for this conjecture based on an analysis of the partial differential equation solved by $F(a,x)$ in the last section. That $1-F$ decays like $q^{2/3}$ can actually be derived from the known Hausdorff dimension (i.e. 4/3) of $\text{SLE}_{8/3}$.  

Our approach rests on the argument of Beffara alluded to above, see Lemma \ref{L:beffara}, and earlier work by Dub\'edat \cite{dubedat:2004}, as well as \cite{RF:rad}, \cite{RF:chord}, where the Loewner equation in multiply connected domains  is discussed and explicit expressions for the change of the conformal parameters under Loewner evolution are given. Using Beffara's argument, it is easy to see that if $D\subset A_q$ is doubly connected, $e^{ix},1\in\partial D$, then
\begin{equation}\label{E:beff2}
P_{A_q,e^{ix}\to1}\{\gamma\subset D\}=\frac{F(a',x')}{F(a,x)}[h'(e^{ix})h'(1)]^{5/8},
\end{equation}
where $h$ is defined in terms of the unique conformal equivalence from $D$ to $A_{q'}$ which keeps 1 fixed, $e^{ix'}$ is the image of $e^{ix}$ under this equivalence, and $a'=\ln q'$. Equation \eqref{E:beff2} is the generalization of \eqref{E:rest4} for $\text{SLE}_{8/3}$.

In \cite{dubedat:2004}, Dub\'edat discusses questions similar to those we discuss here, although he considers $\text{SLE}_6$ and `locality.' Zhan \cite{zhan:ring} constructs $\text{SLE}_2$ in an annulus as the scaling limit of loop-erased random walk, by adapting  the approach taken by Schramm from simply connected domains to doubly connected domains. To do so, he exploits particular properties of the discrete walk. It is also clear from our calculations that $\kappa=2$ is special in that some of the martingales mentioned below have a particularly simple form in this case. However, we will not pursue this here.

Restriction in multiply connected domains has also been discussed in \cite{zhan:thesis}, \cite{dubedat:commutation}, and \cite{lawler:laplacian-b}. In particular, these authors find restriction (local) martingales similar to ours. Due to the greater generality, the expressions these authors find are less explicit and the asymptotics of these (local) martingales  are not discussed. In the case of connectivity two we find here the asymptotics of the restriction martingale, leading to a stochastic functional representation of the intersection probability. We also give a proof that $F$ is smooth enough to apply It\^o's lemma, an issue that, to our knowledge, had not been addressed previously. The  question of smoothness of the intersection probabilities had been raised by John Cardy. While it had been expected that the intersection probability would be given as the solution to a partial differential equation, we are the first to derive it for $\text{SLE}_{8/3}$. A similar equation has been derived in \cite{dubedat:2004} in the percolation case ($\kappa=6$), but the smoothness necessary to apply It\^o's formula in that context is not discussed. Finally, the limiting behavior of the intersection probability as the annulus becomes thinner and thinner is new, though it is clearly related to the estimate obtained in \cite{werner:2005}. It fits with a recent calculation  of Cardy using Coulomb gas methods, \cite{cardy:loop}. 

Going from `locality' to `restriction' in SLE-type calculations involves {\em taking one more derivative}, which leads to expressions which are considerably more expansive. For this reason we begin this paper by changing coordinates from the upper half-plane to a half-strip, where elliptic functions---the indispensable tool of function theory in annuli---have their simplest expression. In Section \ref{S:annulus} we use elliptic functions to describe Loewner evolution in an annulus. In Section \ref{S:mart} we study the `conditional probability martingale' derived from $F$ and use it to show that $F$ has enough smoothness to apply the It\^o formula later in the paper.  In Section \ref{S:asymp} we obtain the asymptotic behavior for $F$ as $q\nearrow1$ and the stochastic representation mentioned above. Finally, in Section \ref{S:pde} we apply It\^o's formula to derive the partial differential equation for $F$.  


\section{Chordal SLE in a half-strip}

Denote $B_t$ a standard 1-dimensional Brownian motion,  $\kappa>0$ 
a constant, and set $W_t=\sqrt{\kappa}B_t$. For $u$ in the upper 
half-plane $\mathbb H$ denote $g_t(u)$ the solution to the chordal 
Loewner equation at time $t$,
\[
\partial_t g_t(u)=\frac{2}{g_t(u)-W_t},\quad g_0(u)=u.
\]
The solution exists up to a time $T_u=\sup \{t:\min_{s\le 
t}|g_s(u)-W_s|>0\}$, and if $K_t=\overline{\{u:T_u\le t\}}$, then 
$g_t$ is the conformal map from $\mathbb H\backslash K_t$ onto 
$\mathbb H$ with hydrodynamic normalization at infinity, 
$\lim_{z\to\infty}g_t(z)-z=0$. The  stochastic process of conformal 
maps $g_t$ is called chordal Schramm-Loewner evolution in $\mathbb H$ 
from $B_0$ to $\infty$ with parameter $\kappa$, see \cite{LSW1}. The 
random growing compact $K_t$ is generated by a curve 
$t\mapsto\gamma_t$ with $\gamma_0=B_0$. If $\kappa\le 4$, then 
$\gamma$ is simple, see \cite{RS:2005}. We will sometimes write 
$\gamma$ for $\gamma[0,\infty)$.

The function
\[
u=\cot (z/2)=i\frac{e^{iz}+1}{e^{iz}-1}
\]
maps the halfstrip $HS\equiv\{z:0\le\Re(z)\le2\pi,\Im(z)<0\}$ onto 
the upper half-plane. We will use $u$ to denote the map as well as 
the variable for the image domain. The sides
\begin{equation}\label{E:sides}
\{iy:y<0\},\quad\{2\pi+iy:y<0\}
\end{equation}
of $HS$
are mapped to the slit $\{iy:y>1\}\subset\mathbb H$ and the real 
interval $(0,2\pi)$ in the $z$-plane corresponds to the real axis in 
the $u$-plane. Furthermore, the point $\infty$ in the (extended) 
$z$-plane corresponds to $i\in\mathbb H$ and the point $\infty$ in 
(the closure of) $\mathbb H$ has the pre-images 
$0,2\pi\in\overline{HS}$. If we identify the sides \eqref{E:sides}, 
i.e. $iy\approx2\pi+iy$, then $u=\cot z/2$ is conformal from $HS$ 
onto $\mathbb H$. In the following we will always assume this 
identification for points in the $z$-plane. The inverse mapping is 
given by
\begin{equation}\label{E:inverse}
z=\frac{1}{i}\ln\frac{u+i}{u-i},
\end{equation}
and we recall the derivatives
\[
u'(z)=-\frac{1}{2}\csc^2(z/2),\quad u''(z)=\frac{1}{2}\csc^2( z/2)\cot (z/2).
\]
We define chordal $\text{SLE}_{\kappa}$ in $HS$ from $x\in(0,2\pi)$ 
to $0$ as the conformal image of  $\text{SLE}_{\kappa}$ in $\mathbb 
H$ from $\cot x/2$ to $\infty$ under the mapping \eqref{E:inverse}.
This definition is natural in light of the characterization of SLE as 
the unique family of measures on non-self-crossing curves which are 
conformally invariant, satisfy a Markovian-type property and a 
certain symmetry condition.

\begin{remark}
It follows from the Riemann mapping theorem that there is a 
one-parameter family of conformal maps from $HS$ onto $\mathbb H$ 
which send $0$ to $\infty$ and $x$ to $\cot x/2$. Choosing a function 
other than $\cot z/2$ from this family would only result in a linear 
time-change for the $SLE$ measures. As we will not be interested in 
{\em when} a particular event occurs but rather {\em if} it occurs 
this is of no concern. In fact, we will change the time parameter 
when it simplifies our calculations.
\end{remark}

If the process $X$ is defined by $X_t=u^{-1}(W_t)$, then
\begin{equation}\label{E:X_t}
	dX_t=-2\sqrt{\kappa}\sin^2 (X_t/2)\ dB_t+2\kappa\sin^4 
(X_t/2)\cot (X_t/2)\ dt.
\end{equation}
Under the random time-change $t\to s$ with $ds=4\sin^4 (X_t/2)\ dt$, we get
\begin{equation}\label{E:X_s}
dX_s=-\sqrt{\kappa}\ dB_s+\frac{\kappa}{2}\cot (X_s/2)\ ds.
\end{equation}
For this new time parameter, let $\tilde{g}_s= u^{-1}\circ g_s\circ 
u$. Then, for each $z\in HS$,
\begin{equation}\label{E:CLEHS}
\partial_s\tilde{g}_s(z)=\Xi_1(\tilde{g}_s(z),X_s),\quad\tilde{g}_0(z)=z,
\end{equation}
with
\[
\Xi_1(z,x)=\frac{2u'(x)^2}{u'(z)[u(z)-u(x)]}=-\frac{\sin^2 
(z/2)}{\sin^4 (x/2)[\cot (z/2)-\cot (x/2)]}.
\]
Note that the vector field $\Xi_1(\cdot,x)$ has a pole with residue 
$2$ at $x$. $\Xi$ is the variation kernel of the Riemann sphere 
expressed in the coordinate $u$, see \cite{SS:1954}. The variation 
kernel is a reciprocal differential (holomorphic vector field) in 
$z$---this explains the $u'$-term in the denominator---and a 
quadratic differential in $x$---which explains the $u'(x)^2$-term in 
the numerator.


\begin{remark}\label{R:stay}
The solution $X_s$ to the SDE \eqref{E:X_s} is a Bessel-like process 
on the interval $(0,2\pi)$. At the boundary points it behaves like 
the 3-dimensional Bessel process, see \cite{lawler:2005}. In 
particular, with probability 1, $X_s$ never leaves $(0,2\pi)$.
\end{remark}


\section{SLE viewed in an annulus}\label{S:annulus}

For a real number $a<0$, $\cot(z/2)$ maps the rectangle \[
R_a\equiv\{0\le\Re(z)\le2\pi,a<\Im(z)<0\}\] onto
$
\mathbb H\backslash C_a$, where $C_a$ denotes the disk
\[
\{u:\left|u-i\frac{1+q^2}{1-q^2}\right|\le\frac{2q}{1-q^2}\},\quad q=e^a.
\]
This doubly connected domain is conformally equivalent to the annulus 
\[A_q\equiv\{q<|z|<1\},\] the image of $R_a$ under the map $z\mapsto 
v=\exp(-iz)$.

For $t>0$, suppose that $\tilde{K}_s\equiv u^{-1}(K_s)\subset R_a$. 
Then $v(\tilde{K}_s)\subset A_q$ and the doubly connected domain 
$A_q\backslash v(\tilde{K}_s)$ is conformally equivalent to a unique 
annulus $A_{q'}$. If $a'=\ln q'$, then $a<a'<0$. Furthermore, there 
is a unique conformal map $\tilde{h}_s:A_q\backslash 
v(\tilde{K}_s)\to A_{q'}$ with $\tilde{h}_s(1)=1$, see 
\cite{ahlfors:1978}. Set $f_s=v^{-1}\circ\tilde{h}_s\circ v$. Then 
$f_s$ maps $R_a\backslash\tilde{K_s}$ onto $R_{a'}$, fixing $0,2\pi$.

To describe the time evolution of $f_s$ we need to use elliptic 
functions. Denote $\zeta$ the Weierstrass $\zeta$-function with 
periods $2\pi$, $2ia$, i.e.
\begin{equation}\label{E:zeta}
\zeta(z)=\zeta(z|a)=\frac{\eta}{\pi}z+\frac{1}{2}\cot 
(z/2)+2\sum_{n=1}^{\infty}\frac{q^{2n}}{1-q^{2n}}\sin nz,
\end{equation}
where
\begin{equation}\label{E:eta}
\eta=\pi\left(\frac{1}{12}-2\sum_{n=1}^{\infty}\frac{n 
q^{2n}}{1-q^{2n}}\right),
\end{equation}
see \cite{hurwitz:1964}. $\zeta$ is regular in the entire $z$-plane 
except for poles with residue 1 at the lattice points $2n\pi+2mia$, 
$n,m\in\mathbb Z$. $\zeta$ is an odd function and $\zeta(\pi)=\eta$.
For each $x\in(0,2\pi)$, $a<0$, define the vector field $\Xi_2(\cdot,x)$ by
\begin{equation}\label{E:Xi_2}
\Xi_2(z,x)=\Xi_2(z,x|a)=2\left[\zeta(z-x)-\frac{\eta}{\pi}z+\zeta(x)\right].
\end{equation}
$\zeta$, $\eta$, and $\Xi_2$ all depend on $a$. We will use $a$ in the notation if any ambiguity as to the particular value of that parameter could arise. 

\begin{proposition}\label{P:Xi_2}
The vector field $\Xi_2(\cdot,x)$ (i) 
is regular except for poles with residue 2 at the points of the 
shifted lattice $\{2n\pi+x+2mia:n,m\in\mathbb Z\}$, (ii) is periodic 
with period $2\pi$ $($i.e. $\Xi_2(z,x)=\Xi_2(z+2\pi,x))$, (iii) 
vanishes at $z=0$, and (iv) has constant imaginary part $+i$, $-i$ on 
the lines $\{\Im(z)=a\}$, $\{\Im(z)=-a\}$, respectively. 
\end{proposition}

\begin{proof}
Property (i) follows immediately from the properties of $\zeta$, and 
(ii), (iii) follow by inspection from \eqref{E:Xi_2}. Next, if 
$\Im(z)=0$, then
\begin{align}
\Im(\Xi_2&(z+ia,x))\notag\\
&=\Im(\cot((z+ia-x)/2))+4\sum_{n=1}^{\infty}\frac{q^{2n}}{1-q^{2n}}\Im(\sin 
n(z+ia-x))\notag\\
&=\frac{1-q^2}{1-2q\cos(z-x)+q^2}-2\sum_{n=1}^{\infty}q^n\cos 
n(z-x)=1,
\end{align}
where the last equality follows from a 
well-known identity for Chebyshev polynomials, see \cite{special:1999}. Similarly, 
$\Im(\Xi_2(z-ia,x))=-1$ if $\Im(z)=0$. 
\end{proof}

For chordal $\text{SLE}_{\kappa}$ in $\mathbb H$ from $0$ to 
$\infty$, and $A<0$, set
\[
	T_A=\inf\{s:K_s\cap C_A\neq\emptyset\}.
\]
If $\kappa\le4$, then $\gamma$ is almost surely a simple curve and thus $K_t=\gamma[0,t]$. In particular, for $\kappa\le 4$, $T_A=\infty$ if and only if $\gamma\cap C_A=\emptyset$.
On $s<T_A$, let $a=a(s)$ be defined as the unique $a$ such that 
\[h_s(\tilde{g}_s(R_A\backslash\tilde{K}_s))= R_{a}.\] Then $a(0)=A$ and $a(s)>a(t)$ if $s>t$ (for an integral expression for $a(s)-a(t)$ see \cite{komatu:1943}). Set 
\[
A^*=\lim_{s\nearrow T_A}a(s).
\] 
Then $A^*\le0$ and $A^*=0$ if and only if $T_A<\infty$. The last statement holds with probability 1 and is a consequence of the fact that a.s. $\gamma_s\to\infty$ as $s\to\infty$. We now change the time parameter from $s$ to $a$ and write $\gamma_a$, $X_a$, $\tilde{g}_a$, and $h_{A,a}$ for $\gamma_{s(a)}$, $X_{s(a)}$, $\tilde{g}_{s(a)}$, and $h_{s(a)}$. We include $A$ in the subscript of $h$ to keep note of the fact that the definition of $h$ depends on $A$ (or rather $R_A$). Then $\gamma[A,a]=\gamma[0,s]$.  

\begin{theorem} For $A\le a<A^*$ we have
\[
	\partial_s a=h_{A,a}'(X_a)^2
\]
and
\begin{equation}\label{E:hdot}
	\partial_a h_{A,a}(z)=\Xi_2(h_{A,a}(z),h_{A,a}(X_a)|a)-\Xi_1(z, 
X_a)\frac{h_{A,a}'(z)}{h_{A,a}'(X_a)^2}.
\end{equation}
\end{theorem}

\begin{proof}
Set $f_{A,a}=h_{A,a}\circ\tilde{g}_a$. Then $f_{A,a}$ is the unique conformal map from $R_A\backslash\gamma[A,a]$ onto $R_a$ with $f_{A,a}(0)=0$. By \cite{komatu:1943},
\begin{equation}\label{E:KL}
\partial_a f_{A,a}(z)=\Xi_2(f_{A,a}(z), Y_{A,a}|a),
\end{equation}
where $Y_{A,a}=h_{A,a}(X_a)$. Note that $Y_{A,a}$ is the image of the tip of the slit $\gamma[A,a]$ under $f_{A,a}$, i.e $Y_{A,a}=\lim_{z\to\gamma_a}f_{A,a}(z)$. Also, it is clear from the mapping properties of $f_{A,a}$ that the left-hand side of \eqref{E:KL} is zero at $z=0$ and has constant imaginary part 1 if $\Im(z)=A$.
Next, by the chain rule
\[
	\partial_a h_{A,a}(z)=\partial_a f_{A,a}(\tilde{g}_a^{-1}(z))
	+(f_{A,a})'(\tilde{g}_a^{-1}(z))\partial_a \tilde{g}_a^{-1}(z).
\]
Since $\partial_a\tilde{g}_a^{-1}(z)=-(\tilde{g}_a^{-1})'(z)(\partial_a\tilde{g}_a)(\tilde{g}_a^{-1}(z))$, we get from \eqref{E:CLEHS}
\[
	\partial_a\tilde{g}_a^{-1}(z)
	=-(\tilde{g}_a^{-1})'(z)\Xi_1(z,X_a)\frac{\partial s}{\partial a}.
\]
Hence
\[
\partial_a h_{A,a}(z)=\Xi_2(h_{A,a}(z),h_{A,a}(X_a)|a)-\Xi_1(z, 
X_a)h_{A,a}'(z)\frac{\partial s}{\partial a},
\]
and
\begin{equation}\label{E:h_s}
\partial_s h_s(z)=\Xi_2(h_{A,a}(z),h_{A,a}(X_a)|a)\frac{\partial a}{\partial s}-\Xi_1(z, 
X_a)h_{A,a}'(z).
\end{equation}
To determine $\partial a/\partial s$ we note that the domains $\tilde{g}_s(R_A)$ change smoothly because $\Xi_1(z,x)$ is smooth away from $x$. The map $h_s$ can be written explicitly in terms of domain functionals, namely the harmonic measures and their conjugates. By Hadamard's formula for the variation of domain functionals under smooth boundary perturbations, see \cite{SS:1954}, it follows that $\partial_s h_s(z)$ extends continuously to the boundary. In particular, the residues of the two terms on the right in \eqref{E:h_s} have to cancel. The residue of the first term is $2(\partial a/\partial s)/h_{A,a}'(X_a)$, the residue of the second $2h'_{A,a}(X_a)$. The theorem now follows.
\end{proof}

We will now draw a number of conclusions from \eqref{E:hdot}. To 
simplify notation, we will indicate differentiation with respect to 
$a$ by $\cdot$, and suppress the subscripts $a, A$ when convenient.

\begin{corollary} \label{C:dh(X)}
On $[A,A^*)$ we have
\begin{equation}\label{E:hdotx}
\dot{h}(X)=2\left[\zeta(h(X))-\frac{\eta}{\pi}h(X)\right]-3\frac{h''(X)}{h'(X)^2}-3\frac{\cot 
(X/2)}{h'(X)}
\end{equation}
and
\begin{align}\label{E:dh(X)}
d(h(X))&=-\sqrt{\kappa}\ 
dB+2\left[\zeta(h(X))-\frac{\eta}{\pi}h(X)\right]\ da\notag\\
&\quad+\frac{\kappa-6}{2}\left[\frac{h''(X)}{ h'(X)^2}+\frac{\cot( 
X/2)}{h'(X)}\right]da.
\end{align}
\end{corollary}

\begin{proof}
Taking the limit $z\to X_a$ in \eqref{E:hdot} gives \eqref{E:hdotx}. The calculation is done by Taylor expansion. By an appropriate version of It\^o's lemma (\cite{revuz.yor:1999}),
\[
d(h(X))=\dot{h}(X)\ da+ h'(X)\ dX+1/2 h''(X)\ dX dX,
\]
where $dX dX$ is the differential of the quadratic variation.
Also, by \eqref{E:X_s},
\begin{equation}\label{E:X_a}
d X_a=-\sqrt{\kappa}\frac{dB}{h'(X)}+\frac{\kappa}{2}\frac{\cot(X/2)}{h'(X)^2}\ da.
\end{equation}
Now \eqref{E:dh(X)} follows from \eqref{E:hdotx}, .
\end{proof}

\begin{remark}
A time change of the results \eqref{E:hdotx} and \eqref{E:dh(X)} had previously been obtained in \cite{dubedat:2004}. 
\end{remark}

Denote $\wp=-\zeta'$ the Weierstrass $\wp$-function,
\[
\wp(z)=\wp(z|a)=-\frac{\eta}{\pi}+\frac{1}{4}\csc^2 
(z/2)-2\sum_{n=1}^{\infty}\frac{n q^{2n}}{1-q^{2n}}\cos nz,
\]
see \cite{hurwitz:1964}.
Then it follows from \eqref{E:hdot} that
\begin{align}\label{E:h'dot}
\dot{h}'(z)&=-2\left[\wp(h(z)-h(X))+\frac{\eta}{\pi}\right]h'(z)\notag\\
&\quad-\frac{h''(z)}{h'(X)^2}\cdot\frac{\sin^3( z/2)}{\sin^3 
(X/2)}\csc\frac{z-X}{2}\notag\\
&\quad+\frac{h'(z)}{h'(x)^2}\frac{\sin^2(z/2)}{\sin^2(X/2)}\left[\frac{1}{2}\csc^2\frac{z-X}{2}-\frac{\cos(z/2)}{\sin(X/2)}\csc\frac{z-X}{2}\right].
\end{align}
In particular,
\begin{equation}\label{E:hprime0}
\dot{h}'(0)=-2\left[\wp(h(X))+\frac{\eta}{\pi}\right]h'(0),
\end{equation}
so that
\begin{equation}\label{E:hprimeexp}
h_{A,a}'(0)=\exp\left(-2\int_A^a\left[\wp(h(X_b))+\frac{\eta}{\pi}\right]db\right).
\end{equation}
Note that $\eta$ in the integrand also depends on $b$, the explicit 
form of the dependence being given in \eqref{E:eta}.

\begin{corollary}
We have
\begin{align}
\dot{h}'(X)&=-\frac{2\eta}{\pi}h'(X)
+\frac{2}{3h'(X)}-\frac{3\cot^2 (X/2)}{2h'(X)}\notag\\
&\quad-3\cot 
(X/2)\frac{h''(X)}{h'(X)^2}+\frac{h''(X)^2}{2h'(X)^3}-\frac{4}{3}\frac{h'''(X)}{h'(X)^2}\notag
\end{align}
and, for real $\alpha$,
\begin{align}
\frac{d(h'(X)^{\alpha})}{\alpha h'(X)^{\alpha}}&=\left[\frac{2}{3 
h'(X)^2}-\frac{3\cot^2 (X/2)}{2h'(X)^2}+\frac{\kappa-6}{2}\cot 
(X/2)\frac{h''(X)}{h'(X)^3}\right]da\notag\\
&\quad+\left[\frac{1+(\alpha-1)\kappa}{2}\frac{h''(X)^2}{h'(X)^4}+
\frac{\kappa-8/3}{2}\cdot\frac{h'''(X)}{h'(X)^3}-\frac{2\eta}{\pi}\right]da\notag\\
&\quad-\sqrt{\kappa}\frac{h''(X)}{h'(X)^2}\ dB.
\end{align}
\end{corollary}

\begin{proof} The first identity follows by taking the limit in \eqref{E:h'dot}, and then the second follows from It\^o's lemma, just as in the proof of Corollary \ref{C:dh(X)}. The calculation is tedious but straightforward and is omitted. 
\end{proof}


\section{Conditional probabilities and restriction martingales}\label{S:mart}

For a simply connected domain $D$ and boundary points $p,q$, we 
define chordal SLE in $D$ from $p$ to $q$ by conformal invariance 
from chordal SLE in $\mathbb H$ from $0$ to $\infty$. This is well 
defined up to a linear time-change. Denote $P_{D,p\to q}$ the law of 
chordal SLE in $D$ from $p$ to $q$, and $\mathbb E_{D,p\to q}$ 
expectation with respect to $P_{D,p\to q}$. Then
\begin{align}\label{E:condition}
&P_{HS,x\to 0}\{\gamma\subset R_A|\gamma[0,s]\}=
P_{\mathbb H,\cot x\to\infty}\{\gamma\cap C_A=\emptyset|\gamma[0,s]\}\notag\\
&\quad=\mathbb E_{\mathbb H, \cot x\to\infty}\left[1\{\gamma[0,s]\cap 
C_A=\emptyset\}1\{\gamma[s,\infty)\cap 
C_A=\emptyset\}|\gamma[0,s]\right]\notag\\
&\quad=1\{s<T_A\}\mathbb E_{\mathbb H, \cot 
x\to\infty}\left[1\{g_s(\gamma[s,\infty))\cap 
g_s(C_A)=\emptyset\}|\gamma[0,s]\right]\notag\\
&\quad=1\{t<T_A\} P_{\mathbb H, W_s\to\infty}\{\gamma\cap 
g_s(C_A)=\emptyset\},
\end{align}
where $W$ is a time changed Brownian motion starting at $\cot x$. We 
note that the last equality follows from basic properties of SLE. Now 
we need a result from \cite{beffara}.

\begin{lemma}{(Beffara)} \label{L:beffara}
Let $\kappa=8/3$.
If $K$ and $K'$ are compact subsets of $\mathbb H$ such that $\mathbb 
H\backslash K$ and $\mathbb H\backslash K'$ are conformally equivalent, then
\[
	P_{\mathbb H,x\to\infty}\{\gamma\cap K=\emptyset\}
	=P_{\mathbb H,\Phi(x)\to\infty}\{\gamma\cap 
K'=\emptyset\}\left[\Phi'(x)\Phi'(\infty)\right]^{5/8},
\]
where $\Phi$ is a conformal map from  $\mathbb 
H\backslash K$ onto $\mathbb 
H\backslash K'$ with $\Phi(\infty)=\infty$ and $\Phi'(\infty)=\lim_{z\to\infty}1/\Phi'(z)$.
\end{lemma}

\begin{theorem} \label{T:prob} If $F(A,x)$ denotes the probability 
that chordal $\rm SLE \sb {8/3}$ in the halfstrip $HS$ from $x$ to 
$0$ stays in the rectangle $R_A$, then
\[
	F(a, h_{A,a}(X_a))\left[\frac{\sin^2 (X_a/2)}{\sin^2 
(h_{A,a}(X_a)/2)}h_{A,a}'(X_a) h_{A,a}'(0)\right]^{5/8}
\]
is a martingale on $[A,A^*)$.
\end{theorem}

\begin{proof} 
It follows from \eqref{E:condition} that $P_{\mathbb H, W_s\to\infty}\{\gamma\cap 
g_s(C_A)=\emptyset\}$ is  a martingale on $s<T_A$. Since
\[
	u\circ h_s\circ u^{-1}(g_s(C_A))=C_a,\quad a=a(s),
\]
it follows from Lemma \ref{L:beffara} that
\begin{align}
&P_{\mathbb H,W_s\to\infty}\{\gamma\cap g_S(C_A)=\emptyset\}\notag\\
&=P_{\mathbb H,u\circ h_s\circ u^{-1}(W_s)\to\infty}\{\gamma\cap C_a=\emptyset\}\notag\\
&\quad\times[(u\circ h_s\circ u^{-1})'(W_s)(u\circ h_s\circ u^{-1})'(\infty)]^{5/8}\notag\\
&=P_{HS,h_a(X_a)\to0}\{\gamma\subset R_a\}[(u\circ h_s\circ u^{-1})'(W_s)(u\circ h_s\circ u^{-1})'(\infty)]^{5/8}.\notag
\end{align}
Next,
\[
(u\circ h\circ u^{-1})'(w)=h'(u^{-1}(w))\sin^2(u^{-1}(w)/2)/\sin^2(h(u^{-1}(w))/2).
\]
If $z=u^{-1}(w)$, then $w\to\infty$ implies $z\to0$. As $\lim_{z\to0}h(z)=0$, we have
\[
\lim_{z\to0}\frac{\sin^2 z/2}{\sin^2(h(z)/2)}=\lim_{z\to0}\left[\frac{\sin z/2}{z}\cdot\frac{z}{h(z)}\cdot\frac{h(z)}{\sin(h(z)/2)}\right]^2=h'(0)^{-2}.
\] 
Since $F(A,x)=P_{HS,x\to0}\{\gamma\subset R_A\}$, the theorem now follows. 
\end{proof}

The martingale in this theorem  is a functional of the Markov process $X_a$ and the non-Markov process $h_{A,a}(X_a)$. Under an appropriate change of measure $h_{A,a}(X_a)$ becomes a Markov process $Y$. This change of measure also introduces a drift to the process in Theorem \ref{T:prob},  and we have to multiply by a factor given by Girsanov's formula to obtain a martingale under this new measure. The new martingale turns out to be a function of $Y$ times an exponential functional of $Y$. Our reason to change measure is that we are able to obtain the asymptotics of this new martingale in Theorem \ref{T:rep}, while it was not clear  to us how to carry out this step for the original martingale in Theorem \ref{T:prob}.   

To change $h_{A,a}(X_a)$ into a Markov process we will first remove the two drift terms in its It\^o decomposition, see \eqref{E:dh(X)}. We carry this out in two steps to better see how the constituent parts fit together. Finally, we perform a third change of measure, which transforms $Y$ from a multiple of a linear Brownian motion to a Bessel-type process on the interval $[0,2\pi]$, a 0-dimensional Legendre process. This last step is natural since it takes the geometry of our setup (i.e. the circle) into account, and, more importantly, leads to a multiplicative stochastic functional in the martingale replacing the martingale from Theorem \ref{T:prob}, whose exponent is an integral with non-singular integrand. 

\begin{proposition} \label{P:m1}If $\kappa=8/3$, $A<0$, and
\[
M_{A,a}=\left[h_{A,a}'(X_a)\frac{\sin^2 (X_a/2)}{\sin^2 
(X_A/2)}\exp\left(\int_A^a\frac{2\eta}{\pi}\ 
db\right)\right]^{5/8},\quad A\le a<A^*,
\]
then $M$ is a martingale with $M_{A,A}=1$ and
\[
dM=-\frac{5}{8}\sqrt{8/3}M\left[\frac{\cot 
(X/2)}{h'(X)}+\frac{h''(X)}{h'(X)^2}\right]dB.
\]
\end{proposition}

\begin{proof}
We have
\begin{align}
&\frac{d\left[h'(X)\sin^2 
(X/2)\right]^{\alpha}}{\alpha\left[h'(X)\sin^2 
(X/2)\right]^{\alpha}}\notag\\
&\quad=-\sqrt{\kappa}\left[\frac{\cot 
(X/2)}{h'(X)}+\frac{h''(X)}{h'(X)^2}\right]dB\notag\\
&\qquad-\frac{2\eta}{\pi}\ 
da+\frac{1+(\alpha-1)\kappa}{2}\cdot\frac{h''(X)^2}{h'(X)^4}\ 
da\notag\\
&\qquad+\left(\frac{8}{3}-\kappa\right)\left[\frac{1}{4h'(X)^2}-\frac{h'''(X)}{2h'(X)^3}\right]da\notag\\
&\qquad+\frac{\kappa(1+2\alpha)-6}{2}\left[\frac{\cot^2 
(X/2)}{2h'(X)^2}+\cot (X/2)\frac{h''(X)}{h'(X)^3}\right]da.
\end{align}
If $\kappa=8/3$ and $\alpha=5/8$ then all drift terms except for the 
first vanish. Since $M$ is also bounded for $a<0$ the proposition 
follows.
\end{proof}

\begin{remark}
If $\kappa>0$ is arbitrary and $\alpha=(6-\kappa)/2\kappa$, then the drift term of $d\left[h'(X)\sin^2 
(X/2)\right]^{\alpha}/\alpha\left[h'(X)\sin^2 
(X/2)\right]^{\alpha}$ reduces to 
\[
-2\eta/\pi\ da+(\kappa-8/3)[\mathcal S h(X)-1/2]/2 h'(X)^2\ da,
\]
where $\mathcal S h=h'''/h'-(3/2) (h''/h')^2$ is the Schwarzian derivative of $h$.
\end{remark}

Denote $P$ the law of the underlying Brownian motion $B$, and denote 
$\mathcal F_a$ the associated filtration after the time-change $t\to 
a$. Define the probability measure $Q$ by
\[
	\frac{dQ}{dP}|_{\mathcal F_a}=M_{A,a}.
\]

\begin{corollary}
Under the measure $Q$,
\[
	F(a, h_{A,a}(X_a))\left[\frac{h_{A,a}'(0)}{\sin^2 (h_{A,a}(X_a)/2)}\exp\left(-\int_A^a\frac{2\eta}{\pi}\ 
db\right)\right]^{5/8}
\]
is a martingale and $Y_{A,a}\equiv h_{A,a}(X_a)$ satisfies
\[
dY=-\sqrt{8/3}\ dB+2\left[\zeta(Y)-\frac{\eta}{\pi}Y\right]da.
\]
\end{corollary}

\begin{proof} The two statements follow from Girsanov's theorem, 
\cite{revuz.yor:1999}, in conjunction with Theorem \ref{T:prob}, 
Proposition \ref{P:m1}, and \eqref{E:dh(X)}.
\end{proof}

Let $\theta(x|a)=\vartheta_1(x/2\pi)$, where $\vartheta_1$ is 
Jacobi's theta function
\[
\theta(x|a)=\vartheta_1(x/2\pi)=-i\sum_{n=-\infty}^{\infty}e^{ix(n+1/2)+a(n+1/2)^2+i\pi 
n}.
\]
Then
\begin{equation}\label{E:heat}
\frac{\partial}{\partial a}\theta(x|a)=-\frac{\partial^2}{\partial 
x^2}\theta(x|a),
\end{equation}
and
\begin{equation}\label{E:thetazeta}
\frac{\partial}{\partial x}\ln\theta(x|a)=\zeta(x)-\frac{\eta}{\pi}x,\quad
\frac{\partial^2}{\partial x^2}\ln\theta(x|a)=-\wp(x)-\frac{\eta}{\pi},
\end{equation}
see \cite{hurwitz:1964}.

We note that if $A^*<0$, then $A^*$ is the first time that $Y$, starting at $Y_A$ at time 
$A<0$, hits $\{0,2\pi\}$. If $A^*=0$, then $Y$ does not hit $\{0,2\pi\}$. 

\begin{proposition} If $\kappa=8/3$, $A<0$, and
\[
N_{A,a}=\left[\vartheta_1\left(\frac{Y_{A,a}}{2\pi}\right)/\vartheta_1\left(\frac{Y_{A,A}}{2\pi}\right)\right]^{-3/4}h_{A,a}'(0)^{1/8},\quad 
A\le a<A^*,
\]
then, under $Q$, $N$ is a martingale with $N_{A,A}=1$ and
\[
dN=\frac{3}{4}\sqrt{8/3} N\left[\zeta(Y)-\frac{\eta}{\pi}Y\right]dB.
\]
\end{proposition}

\begin{proof} Denoting differentiation with respect to the spatial 
variable by a $'$ and using \eqref{E:heat}, we have
\begin{align}
\frac{d\left[\vartheta_1(Y/2\pi)^{\beta}\right]}{\beta\vartheta_1(Y/2\pi)^{\beta}}
&=-\sqrt{\kappa}\frac{\theta'}{\theta}\ dB\notag\\
&\quad+\left[\left(\frac{\kappa}{2}-1\right)\frac{\theta''}{\theta}+\left(2+\frac{\kappa}{2}(\beta-1)\right)\left(\frac{\theta'}{\theta}\right)^2\right]da.
\end{align}
The term in brackets can be rewritten as
\[
\left(1+\frac{\beta\kappa}{2}\right)\frac{\theta''}{\theta}\ 
da+\left[(1-\beta)\frac{\kappa}{2}-2\right](\ln \theta)''\ da.
\]
Thus for $\kappa=8/3$, $\beta=-3/4$,
\begin{align}
d\left[\vartheta_1(Y/2\pi)^{-3/4}\right]&=\frac{3}{4}\sqrt{8/3}\vartheta_1(Y/2\pi)^{-3/4}\frac{\theta'}{\theta}\ 
dB\notag\\
&\quad-\frac{1}{4}\vartheta_1(Y/2\pi)^{-3/4}(\ln \theta)''\ da.
\end{align}
The proposition now follows from \eqref{E:thetazeta} and \eqref{E:hprimeexp}.
\end{proof}

Define the probability measure $R$ for $a<A^*$ by
\[
\frac{dR}{dQ}|_{\mathcal F_a}=N_{A,a}.
\]

\begin{proposition}\label{T:Mart} 
If $Y_{A,a}=h_{A,a}(X_a)$, then under the measure $R$,
\[
F(a,Y_{A,a})\frac{\vartheta_1(Y_{A,a}/2\pi)^{3/4}}{\sin^{5/4}(Y_{A,a}/2)} \exp\left(-\int_A^a\left[\wp(Y_{A,b})+\frac{9\eta}{4\pi}\right]db\right) 
\]
is a martingale for $a<A^*$ and $Y_{A,a}$ satisfies
\[
dY=-\sqrt{8/3}\ dB.
\]
\end{proposition}

\begin{proof}
This is again a consequence of Gisanov's theorem.
\end{proof}

Finally, let 
\[
\tilde{N}_{A,a}=\frac{\sin^{-1/2} (Y_{A,a}/2)}{\sin^{-1/2}(Y_{A,A}/2)}\exp\left[-1/4\int_A^a\csc^2 (Y_{A,b}/2)\ db\right].
\]
It is an easy calculation that---under the measure $R$---$\tilde{N}_{A,a}$ is a martingale on $a<A^*$. If we define the measure $\tilde{R}$ by $d\tilde{R}/dR|\mathcal F_a=\tilde{N}_{A,a}$, then we have the following 

\begin{proposition}\label{T:mart3}
Under the measure $\tilde{R}$ the process $Y_{A,a}$ satisfies
\[
	dY=-\sqrt{8/3}\ dB-2/3\cot Y/2\ da,
\]
and 
\begin{align}\label{E:defmart}
	\mathcal M_{A,a}&\equiv F(a,Y_{A,a})\exp\left[-\int_A^a\left(\wp(Y_{A,b})-\frac{1}{4}\csc^2(Y_{A,b}/2)\right) db\right]\notag\\
	&\quad\times\left[\prod_{n=1}^{\infty}\frac{1-2Q^{2n}+Q^{4n}}{1-2Q^{2n}\cos Y_{A,A}+Q^{4n}}\cdot\frac{1-2q^{2n}\cos Y_{A,a}+q^{4n}}{1-2q^{2n}+q^{4n}}\right]^{3/4}
 \end{align}
is a martingale for $a<A^*$ with $\mathcal M_{A,A}=F(A,Y_{A,A})$. 
\end{proposition}

\begin{proof}
It follows from the infinite product representation of $\vartheta_1$, see \cite{hurwitz:1964}, that 
\begin{equation}\label{E:product}
\frac{\vartheta_1(y/2\pi)}{\sin(y/2)}=q^{1/4}\prod_{n=1}^{\infty}(1-q^{2n})(1-2q^{2n}\cos y+q^{4n}).
\end{equation}
Also,
\[
\exp\left[\sum_{n=1}^{\infty}\int_A^a\frac{2n\tilde{q}^{2n}}{1-\tilde{q}^{2n}}\ db\right]=\prod_{n=1}^{\infty}\frac{1-Q^{2n}}{1-q^{2n}},
\]
and
\begin{align}
&\frac{1-2q^{2n}\cos x+q^{4n}}{1-2Q^{2n}\cos y+Q^{4n}}\notag\\
&=\frac{1-2Q^{2n}+Q^{4n}}{1-2Q^{2n}\cos y+Q^{4n}}\cdot \frac{1-2q^{2n}\cos x+q^{4n}}{1-2q^{2n}+q^{4n}}\left(\frac{1-q^{2n}}{1-Q^{2n}}\right)^2.\notag
\end{align}
Now Girsanov's theorem, Proposition \ref{T:Mart}, and the explicit expression for $\wp$ show that $\mathcal M$ is a martingale. 
\end{proof}

\begin{corollary}\label{C:C12}
For any $A<a<0$, $x\in[0,2\pi]$,
\begin{align}\label{E:1rep}
F(A,x)&=\left(\prod_{n=1}^{\infty}\frac{1-2Q^{2n}+Q^{4n}}{1-2Q^{2n}\cos x+Q^{4n}}\right)^{3/4}\notag\\
&\quad\times\mathbb E[F(a,Y_{A,a})\left(\prod_{n=1}^{\infty}\frac{1-2q^{2n}\cos Y_{A,a}+Q^{4n}}{1-2Q^{2n}+Q^{4n}}\right)^{3/4}\notag\\
&\qquad\times\exp(-\int_A^a\left(\wp(Y)-\frac{1}{4}\csc^2 Y/2\right)db)]
\end{align}
(where $Y_{A,A}=x$), and $F(a,x)$ is $C^{1,2}$ as a function of $a$ and $x$. 
\end{corollary}

\begin{proof}
First, \eqref{E:1rep} is a consequence of \eqref{E:defmart} and the optional sampling theorem. Next, $x\mapsto F(a,x)$ is continuous because the chordal Loewner equation is continuous as a map from the space of continuous paths with the topology of uniform convergence on compacts (the input) to the space of conformal maps with the Caratheodory metric (as output). See \cite{bauer:discrete}, for a discussion. It then follows from the Feynman-Kac formula that the right-hand side of equation \eqref{E:1rep} is $C^{1,2}$, see \cite{KS:1991}.
\end{proof}


\section{Asymptotic behavior of the non-intersection probability}\label{S:asymp}

The stochastic representation of the non-intersection probability
\[
	(a,x)\in[-\infty,0]\times[0,2\pi]\mapsto F(a,x)\equiv P_{\mathbb U,e^{ix}\to1}(\gamma\subset A_q)
\]
we obtain in this section rests on the asymptotics of $F(a,x)$ as $a\nearrow0$. In particular, this probability decays fast enough to control the limiting behavior of the martingale $\mathcal M$ from Proposition \ref{T:mart3}. 

For each $q\in[0,1)$ there exists a unique $L=L(q)\in[0,1)$ such that $A_q$ and $\mathbb U\backslash[-L,L]$ are conformally equivalent. As $q$ increases to 1, $L$ increases to 1 as well. Denote $f$ the conformal equivalence, normalized by $f(1)=1$. For $x\in(0,\pi]$, let $z_1=e^{ix/2}, z_2=e^{-ix/2}$. By symmetry, if $w_{1,2}=f(z_{1,2})$, then $w_2=\bar w_1$.

In what follows we will mean by $h(a)\asymp g(a)$ as $a\nearrow0$, that \[\lim_{a\nearrow0}\log h(a)/\log g(a)=1.
\]

\begin{lemma}\label{L:sn-asymptotics} For $x\in(0,\pi]$, we have
\[
1-L\asymp e^{\frac{\pi^2}{4a}},\text{ and }|f'(z_1)|\asymp|1-f(z_1)|\asymp e^{\frac{\pi}{4a}(\pi-x)}\]
as $a\nearrow0$. 
\end{lemma}

\begin{proof}
From \cite[Chap. VI, Sec. 3]{nehari:1952}, 
\[
f(z)=L\text{ sn}\left(\frac{2iK}{\pi}\log\frac{z}{q}+K;q^4\right),
\]
where $\text{sn}(z)$ is the analytic function for which $\text{sn}'(0)=1$ and which maps the rectangle $\{z:-K<\Re z<K,0<\Im z<i K'\}$ onto the upper half-plane in such a way that $\text{sn}(\pm K)=\pm1$ and $\text{sn}(\pm K+iK')=\pm k^{-1}$. Furthermore, $q^4=\exp(-\pi K'/K)$, and $L=\sqrt{k}$. It is classical that $\text{sn}'(z)=[(1-\text{sn}^2(z))(1-k^2\text{sn}^2(z))]^{1/2}$. Thus
\begin{equation}\label{E:fprime}
f'(z)=(2iK/\pi z)[(L^2-f^2(z))(1-L^2 f^2(z))]^{1/2}.
\end{equation}
Define $h, \tau$ by $q^4=h=e^{i\pi\tau}$, and set $v=\frac{i}{\pi}\log\frac{z_1}{q}+\frac{1}{2}$. Then it follows from \cite[II, 3.]{hurwitz:1964}, and using that texts notation, that 
\[
L=\frac{\theta_2(0|\tau)}{\theta_3(0|\tau)}, \text{ and } f(z)=\frac{\theta_1(v|\tau)}{\theta_0(v|\tau)}.
\]
Using linear transformations of theta functions we may write
\[
\frac{\theta_2(0|\tau)}{\theta_3(0|\tau)}=\frac{\theta_0(0|-\frac{1}{\tau})}{\theta_3(0|-\frac{1}{\tau})},\text{ and } \frac{\theta_1(v|\tau)}{\theta_0(v|\tau)}=i\frac{\theta_1(\frac{v}{\tau}|-\frac{1}{\tau})}{\theta_2(\frac{v}{\tau}|-\frac{1}{\tau})}.
\]
Hence, if $h'=\exp(-i\pi/\tau)$, and using the series representation of $\theta_0$ and $\theta_3$, we get
\[
L=\frac{1+2\sum_{n=1}^{\infty}(-1)^n(h')^{n^2}}{1+2\sum_{n=1}^{\infty}(h')^{n^2}}=1-4h'+O((h')^2),
\]
which is the first statement of the lemma. For the second, we use the infinite product representation of $\theta_1$ and $\theta_2$, giving
\[
i\frac{\theta_1(\frac{v}{\tau}|-\frac{1}{\tau})}{\theta_2(\frac{v}{\tau}|-\frac{1}{\tau})}=\frac{e^{2i\pi v/\tau}-1}{e^{2i\pi v/\tau}+1}\prod_{n=1}^{\infty}\frac{(1-(h')^{2n}e^{2i\pi v/\tau})(1-(h')^{2n}e^{-2i\pi v/\tau})}{(1+(h')^{2n}e^{2i\pi v/\tau})(1+(h')^{2n}e^{-2i\pi v/\tau})}.
\]
Since $\exp(2i\pi v/\tau)=i\exp(-(\pi/4a)(\pi-x))$, the infinite product is $1+O(\exp(\pi^2/(4a)))$, and
\[
\frac{e^{2i\pi v/\tau}-1}{e^{2i\pi v/\tau}+1}=1+2ie^{\frac{\pi}{4a}(\pi-x)}+O(e^{\pi^2/(4a)}),
\]
as $a\nearrow0$. Using equation \eqref{E:fprime}, the lemma now follows. 
\end{proof}

Recall that $z_1=e^{ix/2}$, $w_1=f(z_1)$, and set $u=i(1+w_1)/(1-w_1)$. The following result can be derived from Lemma \ref{L:beffara}. However, we will give a direct argument.

\begin{lemma}\label{L:change} The probability $P_{\mathbb U,e^{ix}\to1}(\gamma\subset A_q)$ is equal to
\[P_{\mathbb H,u\to-u}\left(\gamma\cap i\left[\frac{1-L}{1+L},\frac{1+L}{1-L}\right]=\emptyset\right)\left|\frac{f'(z_1)(1-z_1)}{1-f(z_1)}\right|^{5/4}.
\]
\end{lemma}

\begin{proof}
Denote $B$ a simple curve connecting the inner and outer boundary of $A_q$, so that $B$ is bounded away from $z_1$ and $z_2$. Denote $\phi$ a conformal map from $A_q\backslash B$ onto $\mathbb U$ such that $\phi(z_{1,2})=z_{1,2}$, and $\psi$ a conformal map from $f(A_q\backslash B)$ onto $\mathbb U$ such that $\psi(w_{1,2})=w_{1,2}$. Then, by conformal restriction, 
\begin{align}\label{E:psis}
P_{\mathbb U,z_1\to z_2}(\gamma\subset A_q\backslash B)&=|\phi'(z_1)\phi'(z_2)|^{5/8},\notag\\
P_{\mathbb U,w_1\to w_2}(\gamma\subset f(A_q\backslash B))&=|\psi'(w_2)\psi'(w_2)|^{5/8}. 
\end{align}
Since $T\equiv\phi\circ f\circ \psi^{-1}$ maps $\mathbb U$ onto $\mathbb U$ and sends $w_{1,2}$ to $z_{1,2}$, there is a  pair $w_0, z_0\in\partial\mathbb U$ such that $T$ is the linear transformation given by
\[
\frac{T(w)-w_1}{T(w)-w_2}\cdot\frac{w_0-w_2}{w_0-w_1}=\frac{z-z_1}{z-z_2}\cdot\frac{z_0-z_2}{z_0-z_1}.
\]
A calculation gives
\[
T'(w_1)T'(w_2)=\left(\frac{z_1-z_2}{w_1-w_2}\right)^2,
\]
which together with $|f'(z_1)|=|f'(z_2)|$ implies
\begin{equation}\label{E:part}
P_{\mathbb U,z_1\to z_2}(\gamma\subset A_q\backslash B) =P_{\mathbb U,w_1\to w_2}(\gamma\subset f(A_q\backslash B))\left|\frac{f'(z_1)(z_1-z_2)}{w_1-w_2}\right|^{5/4}.
\end{equation}
By an inclusion/exclusion argument, equation \eqref{E:part} also holds if $A_q\backslash B$ is replaced by $A_q$. Finally, by conformal invariance, 
\[
P_{\mathbb U,w_1\to w_2}(\gamma\subset f(A_q))=P_{\mathbb H,u\to-u}\left(\gamma\cap i\left[\frac{1-L}{1+L},\frac{1+L}{1-L}\right]=\emptyset\right).
\]
\end{proof}

Note that because $x\in(0,\pi]$ we have $\arg z_1,\arg w_1\in(0,\pi/2]$ and so $u\le-1.$ We will use the following lower and upper bounds:
\begin{align}\label{E:lower}
&P_{\mathbb H,u\to-u}\left(\gamma\cap i\left[\frac{1-L}{1+L},\frac{1+L}{1-L}\right]=\emptyset\right)\notag\\
&\ge P_{\mathbb H,u\to-u}(\gamma\cap i(0,\frac{1+L}{1-L}]=\emptyset)+P_{\mathbb H,u\to-u}(\gamma\cap i[\frac{1-L}{1+L},\infty)=\emptyset)\notag\\
&=P_{\mathbb H,u\to-u}(\gamma\cap i(0,\frac{1+L}{1-L}]=\emptyset)+P_{\mathbb H,\frac{1}{u}\to-\frac{1}{u}}(\gamma\cap i(0,\frac{1+L}{1-L}]=\emptyset),
\end{align}
and
\begin{align}\label{E:upper}
&P_{\mathbb H,u\to-u}\left(\gamma\cap i\left[\frac{1-L}{1+L},\frac{1+L}{1-L}\right]=\emptyset\right)\notag\\
&\le P_{\mathbb H,u\to-u}(\gamma\cap i(0,\frac{1+L}{1-L}]=\emptyset)+P_{\mathbb H,\frac{1}{u}\to-\frac{1}{u}}(\gamma\cap i[\frac{1+L}{1-L},\infty)=\emptyset)\notag\\
&\quad+P_{\mathbb H,u\to-u}(\gamma\cap i(0,\frac{1-L}{1+L})\neq\emptyset,\gamma\cap i(\frac{1+L}{1-L},\infty)\neq\emptyset).
\end{align}

For $c\in\mathbb R, d>0$, set 
\[
g_{c,d}(z)=\frac{|c|}{\sqrt{c^2+d^2}}\sqrt{z^2+d^2}.
\]
Then $g_{c,d}$ maps $\mathbb H\backslash i(0,d]$ conformally onto $\mathbb H$ such that $g_{c,d}(\pm c)=\pm c$. Furthermore,
\[
|g'_{c,d}(c)g'_{c,d}(-c)|=\frac{c^4}{(c^2+d^2)^2},
\]
and so by conformal restriction
\begin{equation}\label{E:gcd}
P_{\mathbb H,c\to -c}(\gamma\cap i(0,d]=\emptyset)=[c^2/(c^2+d^2)]^{5/4}.
\end{equation}

\begin{corollary}\label{C:gcd}
We have
\[
P_{\mathbb H,u\to-u}(\gamma\cap i(0,\frac{1+L}{1-L}]=\emptyset)+P_{\mathbb H,\frac{1}{u}\to-\frac{1}{u}}(\gamma\cap i(0,\frac{1+L}{1-L}]=\emptyset)\asymp e^{\frac{5\pi x}{8a}}
\]
as $a\nearrow0$. 
\end{corollary}

\begin{proof}
By \eqref{E:gcd},
\[
P_{\mathbb H,u\to-u}(\gamma\cap i(0,\frac{1+L}{1-L}]=\emptyset)
=\left(\frac{u(1-L)}{1+L}\right)^{5/2}\left(1+\frac{(u(1-L))^2}{(1+L)^2}\right)^{-5/4},
\]
and from Lemma \ref{L:sn-asymptotics}
\[
\left(\frac{u(1-L)}{1+L}\right)^{5/2}\left(1+\left(\frac{u(1-L)}{1+L}\right)^2\right)^{-5/4}\asymp e^{\frac{5\pi x}{8a}}.
\]
Similarly,
\[
P_{\mathbb H,\frac{1}{u}\to-\frac{1}{u}}(\gamma\cap i(0,\frac{1+L}{1-L}]=\emptyset)\asymp e^{\frac{5\pi^2}{8a}+\frac{5\pi}{8a}(\pi-x)},
\]
so that this term is negligible compared to the first if $0<x<\pi$, and of the same order if $x=\pi$.
\end{proof}

\begin{lemma}\label{L:upper-lower} We have
\[
P_{\mathbb H,u\to-u}\left(\gamma\cap i(0,\frac{1-L}{1+L})\neq\emptyset,\gamma\cap i(\frac{1+L}{1-L},\infty)\neq\emptyset\right)\asymp e^{\pi^2/a},
\]
as $a\nearrow0$.
\end{lemma}

\begin{proof} First,
\begin{align}\label{E:first}
&P_{\mathbb H,u\to-u}\left(\gamma\cap i(0,\frac{1-L}{1+L})\neq\emptyset,\gamma\cap i(\frac{1+L}{1-L},\infty)\neq\emptyset\right)\notag\\
&=P_{\mathbb H,u\to-u}\left(\gamma\cap i(0,\frac{1-L}{1+L})\neq\emptyset\right)+P_{\mathbb H,\frac{1}{u}\to-\frac{1}{u}}\left(\gamma\cap i(0,\frac{1-L}{1+L})\neq\emptyset\right)\notag\\
&\quad-P_{\mathbb H,u\to-u}\left(\gamma\cap i((0,\frac{1-L}{1+L})\cup(\frac{1+L}{1-L},\infty))\neq\emptyset\right).
\end{align}
The last probability on the right equals
\[
P_{\mathbb U,w_1\to w_2}(\gamma\cap((-1,-L]\cup[L,1))\neq\emptyset).
\]
To calculate this probability, note that 
\[
g_L(w)\equiv\frac{1+w^2-\sqrt{(1+w^2)^2-4p^2w^2}}{2pw}
\]
maps $\mathbb U\backslash((-1,-L]\cup[L,1))$ onto $\mathbb U$ if $2p=(L+1/L)$, see \cite[Chapter 3]{ivanov}. Here, the square root is chosen so that $g_L(i)=i$. Setting $w=e^{i\varphi}$, this can be written
\begin{equation}
g_L(w)=\begin{cases}\frac{1}{p}\cos\varphi+i\sqrt{1-\frac{1}{p^2}\cos^2\varphi}, &\text{if $\varphi\in(0,\pi/2]$;}\\
\frac{1}{p}\cos\varphi-i\sqrt{1-\frac{1}{p^2}\cos^2\varphi}, &\text{if $\varphi\in[-\pi/2,0)$.}
\end{cases}
\end{equation}
Then
\[
g'_L(w)g'_L(\bar w)=-\frac{\sin^2\varphi}{p^2-1+\sin^2\varphi}.
\]
Denote $T$ a (fractional) linear transformation from $\mathbb U$ onto $\mathbb U$ sending $g_L(w_{1,2})$ onto $w_{1,2}$. Then, as in the proof of Lemma \ref{L:change}, 
\[
T'(g_L(w_1))T'(g_L(w_2))=\frac{\sin^2\varphi}{1-\frac{1}{p^2}\cos^2\varphi},
\]
where now $\varphi=\arg w_1$. Thus, by conformal restriction,
\begin{equation}\label{E:p1}
P_{\mathbb U,w_1\to w_2}(\gamma\cap((-1,-L]\cup[L,1))\neq\emptyset)=1-\left[\frac{p\sin^2\varphi}{p^2-1+\sin^2\varphi}\right]^{5/4}.
\end{equation}
Finally, from the definition of $u$ and $\varphi$ in terms of $w_1$, it follows that $u=-\cot(\varphi/2)$ and so $4/\sin^2\varphi=(u+1/u)^2$. A calculation now gives
\begin{align}\label{E:p2}
&\frac{p^2-1+\sin^2\varphi}{p\sin^2\varphi}\notag\\
&=1+\left(\frac{1-L}{1+L}\right)^2(u^2+\frac{1}{u^2})+\frac{(1-L)^4}{8(L+L^3)}\left[2+\left(\frac{1-L}{1+L}\right)^2(u^2+\frac{1}{u^2})\right].
\end{align}
On the other hand, \eqref{E:gcd} implies
\begin{equation}\label{E:hit1}
P_{\mathbb H,u\to-u}\left(\gamma\cap i(0,\frac{1-L}{1+L})\neq\emptyset\right)=1-\left(1+\left(\frac{1-L}{1+L}\right)^2\frac{1}{u^2}\right)^{-5/4}
\end{equation}
and
\begin{equation}\label{E:hit2}
P_{\mathbb H,\frac{1}{u}\to-\frac{1}{u}}\left(\gamma\cap i(0,\frac{1-L}{1+L})\neq\emptyset\right)=1-\left(1+\left(\frac{1-L}{1+L}\right)^2 u^2\right)^{-5/4}.
\end{equation}
Combining \eqref{E:hit1}, \eqref{E:hit2}, \eqref{E:p1}, and \eqref{E:first}, we get
\begin{align}\label{E:ende}
&P_{\mathbb H,u\to-u}\left(\gamma\cap i(0,\frac{1-L}{1+L})\neq\emptyset,\gamma\cap i(\frac{1+L}{1-L},\infty)\neq\emptyset\right)\notag\\
&=1-\left(1+\left(\frac{1-L}{1+L}\right)^2\frac{1}{u^2}\right)^{-5/4}
+1-\left(1+\left(\frac{1-L}{1+L}\right)^2 u^2\right)^{-5/4}\notag\\
&\quad-1+\left(\frac{p^2-1+\sin^2\varphi}{p\sin^2\varphi}\right)^{-5/4}.
\end{align}
Using \eqref{E:p2}, straightforward expansion of the right hand side of \eqref{E:ende} shows it to be equal to 
\[
\frac{5}{256}(1-L)^4+\frac{5}{128}(1-L)^5+(1-L)^4O(u^2(1-L)^2).
\]
\end{proof}

From the upper and lower bounds \eqref{E:upper}, \eqref{E:lower}, Corollary \ref{C:gcd} and Lemma \ref{L:upper-lower} we get

\begin{theorem}\label{T:asymp0} For every $x\in(0,\pi]$ we have
\begin{equation}\label{E:asymp0}
F(a,x)\asymp\exp\left(\frac{5\pi x}{8a}\right)
\end{equation}
as $a\nearrow0$.
\end{theorem}

We now combine the previous result and Proposition \ref{T:mart3} to obtain a stochastic representation of $F(A,x)$.

\begin{theorem}\label{T:rep}
Under the measure $\tilde{R}$ we have \[\sup_{a<A^*}|\mathcal M_{A,a}|<\infty.\] Furthermore, if $A^*=0$, then $\lim_{a\nearrow A^*}\mathcal M_{A,a}=0$, while if $A^*<0$ and $Q^*=e^{A^*}$, then
\begin{align}\label{E:limmart}
\lim_{a\nearrow A^*}\mathcal M_{A,a}&=\left[\prod_{n=1}^{\infty}\frac{1-2Q^{2n}+Q^{4n}}{1-2Q^{2n}\cos Y_{A,A}+Q^{4n}}\right]^{3/4}\notag\\
&\qquad\times\exp\left[-\int_A^{A^*}\left(\wp(Y_{A,b})-\frac{1}{4}\csc^2(Y_{A,b}/2)\right) db\right].
\end{align}
Finally, if $x=Y_{A,A}$, then  
\begin{align}\label{E:rep}
F(A,x)&=\left[\prod_{n=1}^{\infty}\frac{1-2Q^{2n}+Q^{4n}}{1-2Q^{2n}\cos x+Q^{4n}}\right]^{3/4}\notag\\
&\quad\times\mathbb E^{\tilde{R}}\left[\exp\left[-\int_A^{A^*}(\wp(Y_{A,b})-\frac{1}{4}\csc^2(Y_{A,b}/2)) db\right], A^*<0\right].
\end{align}
\end{theorem}

\begin{proof}
That $\mathcal M$ is a bounded martingale follows from the limiting behavior as $a\nearrow A^*$, which we now establish. If $A^*<0$, then $Y_{A,A^*}=0$ and \eqref{E:limmart} follows directly from \eqref{E:defmart}. On the other hand, if $A^*=0$, then $Y_{A,A^*}\neq0$ a.s. and it follows from Theorem \ref{T:asymp0} that $F(a,Y_{A,a})$ decays like $\exp(-cx/(1-q))$ with $c=5\pi/8$, and $x=\min\{Y_{A,A},2\pi-Y_{A,A}\}$. We will now show that  
\begin{align}\label{E:estprod}
&\left[\prod_{n=1}^{\infty}\frac{1-2q^{2n}\cos x+q^{4n}}{1-2q^{2n}+q^{4n}}\right]^{3/4}\notag\\
&\qquad\le\exp\left[\frac{1}{1-q}\left(\frac{\pi^2}{8}-\frac{3}{8}[\text{Li}_2(e^{ix})+\text{Li}_2(e^{-ix})]\right)\right],
\end{align}
where $\text{Li}_2$ denotes the dilogarithm. 
Set $x_n=1-q^{2n}$, $n\ge0$. Then $x_n-x_{n-1}=(1-x_n)(1-q^2)/q^2$, and by simple integral comparison,
\begin{align}
\sum_{n=1}^{\infty}\ln(1-q^{2n})&=\frac{q^2}{1-q^2}\sum_{n=1}^{\infty}\frac{\ln x_n}{1-x_n}(x_n-x_{n-1})\notag\\
&\ge\frac{q^2}{1-q^2}\int_0^1\frac{\ln x}{1-x}\ dx=-\frac{1}{1-q}\cdot\frac{\pi^2q^2}{6(1+q)}.\notag
\end{align}
Thus
\begin{equation}\label{E:denominator}
-\frac{3}{2}\sum_{n=1}^{\infty}\ln(1-q^{2n})\le\frac{1}{1-q}\cdot\frac{\pi^2}{8}.
\end{equation}
Similarly, if we set $y_n=-q^{2n}$, $n\ge0$, then $y_{n}-y_{n-1}=-y_n(1-q^2)/q^2$ for $n\ge1$, and so
\begin{align}\label{E:numerator}
&\sum_{n=1}^{\infty}\ln(1-2q^{2n}\cos x+q^{4n})\notag\\
&=-\frac{q^2}{1-q^2}\sum_{n=1}^{\infty}\frac{\ln (1+2y_n\cos x+y_n^2)}{y_n}(y_{n}-y_{n-1})\notag\\
&\le\frac{q^2}{1-q^2}\int_0^1\frac{\ln(1-2y\cos x+y^2)}{y} dy \notag\\
&=\frac{q^2}{1-q^2}\left(-\text{Li}_2(e^{ix})-\text{Li}_2(e^{-ix})\right).
\end{align}
Now, \eqref{E:denominator} and \eqref{E:numerator} imply \eqref{E:estprod}. It is elementary that 
\[
5\pi x\ge\pi^2-3[\text{Li}_2(e^{ix})+\text{Li}_2(e^{-ix})]
\]
for $x\in[0,\pi]$, with equality holding for $x=0$. 
Thus $\mathcal M$ is a bounded martingale and \eqref{E:rep} follows from the optional sampling theorem.
\end{proof}

\begin{remark}
Under $\tilde{R}$, $Y$ is a Legendre process whose boundary behavior is that of a 0-dimensional Bessel process, i.e. 0 and $2\pi$ are absorbing, see \cite{revuz.yor:1999}. It can also be interpreted as the driving function of a radial $\text{SLE}(\kappa, \rho)$.
By \eqref{E:product},
\[
\prod_{n=1}^{\infty}\frac{1-2Q^{2n}+Q^{4n}}{1-2Q^{2n}\cos x+Q^{4n}}\]
is the quotient of $y\mapsto\vartheta_1(y/2\pi)/\sin(y/2)$ evaluated at $y=0$ and at $y=x$. Also, 
\begin{align}
\exp&\left[-\int_A^{A^*}\left(\wp(Y_{A,a})-\frac{1}{4}\csc^2(Y_{A,a}/2)\right) da\right]\notag\\
&=\left(\frac{Q^*}{Q}\right)^{1/12}\exp\left[-\int_A^{A^*}\frac{2nq^{2n}}{1-q^{2n}}\left(1-\cos n Y_{A,a}\right) da\right].
\end{align}
\end{remark}

\begin{remark} Obviously, $F(a,x)=\mathbb E^P[1,A^*<0]$, where $P$ is the original SLE-measure under which $Y=h(X)$ is the non-Markov process satisfying equation \eqref{E:dh(X)}. Thus the price we incur   for switching to a Markov process representation is an exponential functional. We note that this exponential functional can be given an interpretation using the Brownian loop soup. 
\end{remark}


\section{The Partial differential equation}\label{S:pde}

It follows from Corollary \ref{C:C12} that $F(a,x)$ is smooth enough in $(a,x)$ to apply It\^o's formula, and we have

\begin{theorem}\label{T:main}
If $G(a,x)=F(a,x)\vartheta_1(x/2\pi)^{3/4}\sin^{-5/4}(x/2)$, then
\begin{equation}\label{E:pde-short}
-\partial_a G=\frac{4}{3}G''-\left(\wp(x)+\frac{9\eta}{4\pi}\right)G.
\end{equation}
Furthermore, $F(a,x)$ is the unique solution to the evolution equation
\begin{align}\label{E:pde-full}
-\partial_a 
F&=\frac{4}{3}F''+\left[2\zeta(x)-\frac{2\eta}{\pi}x-\frac{5}{3}\cot\frac{x}{2}\right]F'\notag\\
&\quad+\left[\frac{15}{16}\csc^2\frac{x}{2}-\frac{5}{4}\left(\cot\frac{x}{2}\left[\zeta(x)-\frac{\eta}{\pi}x\right]+\wp(x)+\frac{2\eta}{\pi}+\frac{5}{12}\right)\right]F,
\end{align}
for $(a,x)\in(-\infty,0)\times(0,2\pi)$, and
with initial condition \[\lim_{a\searrow-\infty}F(a,x)=1,\]
 and boundary condition \[F(a,0)=F(a,2\pi)=1.\] Finally, the solution $F$ is symmetric, $F(a,x)=F(a,2\pi-x)$. 
\end{theorem}

\begin{proof}
The partial differential equation for $G$ is a consequence of Theorem 
\ref{T:Mart} and It\^o's lemma. The equation for $F$ follows from the 
equation for $G$. Finally, that $F(a,0)=1$ is clear and it is also known, for example by considering the Hausdorff dimension of the $\text{SLE}_{8/3}$ curve, that $\lim_{a\to-\infty}F(a,x)=1$.
\end{proof}

We now briefly discuss the case $q\searrow0$. As we could not find stronger convergence results for PDEs such as \eqref{E:pde-full} in the literature we can only establish the rate in a weak sense, see Remark \ref{R:galerkin}. 

Using the formulas for $\zeta$, $\eta$, and $\wp$, we can write \eqref{E:pde-full} as
\begin{align}\label{E:pde-sum}
-\partial_a 
F&=\frac{4}{3}F''+\left[-\frac{2}{3}\cot(x/2)+4\sum_{n=1}^{\infty}\frac{q^{2n}}{1-q^{2n}}\sin nx\right]F'\notag\\
&\quad+\frac{5}{2}\sum_{n=2}^{\infty}\frac{q^{2n}}{1-q^{2n}}\left[n(1+\cos nx)-\cot(x/2)\sin nx\right]\cdot F.
\end{align}
In particular, the coefficient of the zeroth-order term is nonsingular in $x$ and vanishes for $x=0$. We note also that the summation in the zeroth-order term begins with $n=2$ because $(1+\cos x)/\sin x=\cot x/2$.  

To guess the behavior of $F$ as $q\searrow0$ we consider the PDE obtained by setting $q=0$ in \eqref{E:pde-sum}, 
\begin{equation}\label{E:pde-asymp}
-\partial_a H=4/3\ H''-2/3\ \cot x/2 \ H'.
\end{equation}
Then \eqref{E:pde-sum} is a perturbation of \eqref{E:pde-asymp} if $q$ is small.
If we replace $H$ by $1-H$, then $1-H$ satisfies the same equation. We consider the mixed initial-boundary value problem for \eqref{E:pde-asymp} where  
\begin{equation}\label{E:initial}
\lim_{a\to-\infty}H(a,x)=0,\text{ for }x\in(0,2\pi),\text{ and } H(a,0)=0, \text{ for }a\in(-\infty,0).
\end{equation}
The solution should describe the asymptotic behavior of $P_{\mathbb H,\cot x\to\infty}\{\gamma\cap C_a\neq\emptyset\}$ as $a\to-\infty$.

\begin{proposition}\label{P:asymp}
The solutions to the mixed initial-boundary value problem \eqref{E:pde-asymp},\eqref{E:initial} are given by
\[
	H(a,x)=c q^{2/3}\sin^2 x/2,
\]
for an arbitrary positive constant $c$.
\end{proposition}

\begin{proof} This follows easily from separation of variables. 
\end{proof}

\begin{remark}
The exponent $2/3$ is as expected. It is a special case of the ``first moment estimate'' given in \cite{beffara:hausdorff}, where it is shown that the Hausdorff dimension of $\text{SLE}_{8/3}$ is 4/3.
\end{remark}
  
It is clear from the form of the equation \eqref{E:pde-asymp} and the initial-boundary value conditions that multiplication of a solution by a constant gives another solution. For the full equation \eqref{E:pde-sum} this is not the case. The corresponding equation for $1-F$ has the same initial and boundary value conditions as \eqref{E:initial} but the equation is no longer homogeneous.    

\begin{remark}\label{R:galerkin}
The Galerkin approximation, see \cite{evans:pde}, for \eqref{E:pde-sum} (or rather for the inhomogeneous equation satisfied by $1-F$), using the orthonormal system $(1/\sqrt{\pi})\sin((2k-1)x/2)$, $k=1,2,\dots$, gives as first approximation to $1-F$
\[
	{\pi}^{-1/2}q^{2/3}(1-q^2)^{1/2}\prod_{n=2}^{\infty}(1-q^{2n})^{5/4}\sin(x/2).
\]
It is a weak solution of the equation for $1-F$ when testing against the 1-dimensional space spanned by $w_1$. For larger subspaces the systems of ODEs the Galerkin approximation gives rise to, did not appear tractable to us.
\end{remark}

\end{document}